\documentclass[11pt,a4paper]{article}
\usepackage{amssymb}
\usepackage{amsmath}
\usepackage{harvard}
\usepackage{cite}

\input{tcilatex}
\begin{document}

\title{Regular circle actions on 2-connected 7-manifolds}
\author{Yi Jiang\thanks{%
The author's research is partially supported by 973 Program 2011CB302400 and
NSFC 11131008.} \\
Institute of Mathematics, Chinese Academy of Sciences, \\
Beijing 100190, China}
\maketitle

\begin{abstract}
We determine the homeomorphism (resp. diffeomorphism) types of those
2-connected 7-manifolds (resp. smooth 2-connected 7-manifolds) that admit
regular circle actions (resp. smooth regular circle actions).

\begin{description}
\item[2010 Mathematics subject classification:] 55R15 (55R40)

\item[Key Words and phrases:] group actions, Kirby--Siebenmann invariant, $%
\mu -$invariant, characteristic classes

\item[Email address:] jiangyi@amss.ac.cn
\end{description}
\end{abstract}

\section{Introduction}

In this paper all manifolds under consideration are closed, oriented and
topological, unless otherwise stated. Moreover, all homeomorphisms and
diffeomorphisms are to be orientation preserving. Given a positive integer $%
n $ let $S^{n}$ (resp. $D^{n+1}$) be the diffeomorphism type of the unit $n$%
--sphere (resp. the unit $(n+1)$--disk) in $(n+1)$--dimensional Euclidean
space $\mathbb{%
\mathbb{R}
}^{n+1}$.

\bigskip

\noindent \textbf{Definition 1.1.} A circle action $S^{1}\times M\rightarrow
M$ on a manifold $M$ is called \textsl{regular} if this action is free and
the orbit space $N:=M/S^{1}$ (with quotient topology) is a manifold.

Similarly, a smooth circle action $S^{1}\times M\rightarrow M$ on a smooth
manifold $M$ is called \textsl{regular }if this action is free (see \cite[%
p.38, Proposition 5.2]{T}).

\bigskip

For a given manifold $M$ one can ask

\noindent \textbf{Problem 1.2. }Does $M$ admit a regular circle action?

\bigskip

Solutions to Problem 1.2 can have direct implications in contact topology.
For example, the Boothby-Wang theorem implies that the existence of a smooth
regular circle action on a smooth manifold $M$ is a necessary condition to
the existence of a regular contact form on $M$ (see \cite[p.341]{G}).

Problem 1.2 has been solved for all $1$--connected $5$--manifolds by Duan
and Liang\cite{Du}. In particular, it was shown that all $1$--connected $4$%
--manifolds with second Betti number $r$ can be realized as the orbit spaces
of some regular circle actions on the single $5$--manifold $%
\#_{r-1}S^{2}\times S^{3}$, the connected sums of $r$--$1$ copies of the
product $S^{2}\times S^{3}$. In this paper we study Problem 1.2 for the $2$%
--connected $7$--manifolds.

Our main result is stated in terms of a family $M_{l,k}^{c},c\in \left\{
0,1\right\} ,l,k\in 
\mathbb{Z}
$ of $2$--connected $7$--manifolds. The manifolds $M_{l,k}^{0}$ are the
total spaces of the $S^{3}$--bundles $\pi _{M}:M_{l,k}^{0}\rightarrow S^{4}$
with characteristic map $[f_{l,k}]\in \pi _{3}(SO_{4})$ defined by

\begin{quote}
$f_{l,k}(u)v=u^{l+k}vu^{-l}$, $v\in \mathbb{%
\mathbb{R}
}^{4}$, $u\in S^{3}$
\end{quote}

\noindent where the space $\mathbb{%
\mathbb{R}
}^{4}$ and the sphere $S^{3}$ are identified with the algebra of quaternions
and the space of unit quaternions respectively, and where quaternion
multiplication is understood on the right hand side of the formula. Complete
classification on the manifolds $M_{l,k}^{0}$ has been obtained by Milnor%
\cite{M}, Crowley and Escher\cite{C}.

The manifold $M_{l,k}^{1}$ is an analogue of the manifold $M_{l,k}^{0}$ in
the non--smooth category when $k\equiv 0\func{mod}2$. Write $\mathcal{S}%
^{TOP}(M_{l,k}^{0})$ for the set of all equivalence classes $[M,h]$ of the
pairs $(M,h)$ with $h:M\rightarrow M_{l,k}^{0}$ a homotopy equivalence of $7$%
--manifolds. Two such pairs $(M,h)$ and $(M^{^{\prime }},h^{^{\prime }})$
are called \textsl{equivalent} if there is a homeomorphism $f:M\rightarrow
M^{^{\prime }}$ such that $h$ is homotopic to $h^{^{\prime }}\circ f$ (see 
\cite[Chapter 2]{Ma}). Adapting the arguments of \cite[Section 5]{C} from
the PL case to the TOP case we have the composition of isomorphisms

\begin{center}
$\mathcal{S}^{TOP}(M_{l,k}^{0})\overset{\eta }{\rightarrow }%
[M_{l,k}^{0},G/TOP]\overset{d}{\rightarrow }H^{4}(M_{l,k}^{0};\pi
_{4}(G/TOP))\cong H^{4}(M_{l,k}^{0})$,
\end{center}

\noindent where the space $G$ (resp. $TOP$) is the direct limit of the set
of self homotopy equivalences of $S^{n-1}$ (resp. the topological monoid of
origin-preserving homeomorphisms of $%
\mathbb{R}
^{n}$), $\eta $ is the one to one correspondence in the surgery exact
sequence (see \cite[p.40-44]{Ma}) and where the isomorphism $d$ is induced
by the primary obstruction to null-homotopy. Write $\left[ M_{l,k}^{1},h_{M}%
\right] $ for the element $(d\circ \eta )^{-1}(\pi _{M}^{\ast }(\iota ))\in 
\mathcal{S}^{TOP}(M_{l,k}^{0})$ where $\iota \in H^{4}(S^{4})$ is the
generator as in \cite{C} and $H^{4}(M_{l,k}^{0})\cong 
\mathbb{Z}
_{k}$ is generated by $\pi _{M}^{\ast }(\iota )$. Clearly, the 7-manifold $%
M_{l,k}^{1}$ is $2$--connected and unique up to homeomorphism.

The group $\Gamma _{7}$ of exotic $7$--spheres is cyclic of order $28$ with
generator $M_{1,1}^{0}$ \cite[Section 6]{E}. Let $\Sigma
_{r}:=rM_{1,1}^{0}\in \Gamma _{7}$, $r\in 
\mathbb{Z}
$. Our main result is stated below, where $\mathbb{%
\mathbb{N}
}$ is the set of all nonnegative integers.

\bigskip

\noindent \textbf{Theorem 1.3. }All homeomorphism classes of the $2$%
--connected $7$--manifolds that admit regular circle actions are represented
by

\begin{center}
$\#_{2r}S^{3}\times S^{4}\#M_{6m,(1+c)k}^{c}$, $c\in \left\{ 0,1\right\}
,r\in \mathbb{%
\mathbb{N}
}$ and $m,k\in \mathbb{%
\mathbb{Z}
}$.
\end{center}

All diffeomorphism classes of the smooth $2$--connected $7$--manifolds that
admit smooth regular circle actions are represented by

\begin{center}
$\#_{2r}S^{3}\times S^{4}\#M_{6(a+1)m,(a+1)k}^{0}\#\Sigma _{(1-a)m}$, $a\in
\left\{ 0,1\right\} $, $r\in \mathbb{%
\mathbb{N}
}$, $m,k\in \mathbb{%
\mathbb{Z}
}$.
\end{center}

In the course to establish Theorem 1.3 we obtain also a classification of
the $6$--manifolds that can appear as the orbit spaces of some regular
circle actions on $2$--connected $7$--manifolds, see Lemmas 2.1 and 2.2 in
Section 2. In addition, Theorem 1.3 has some direct consequences which are
discussed in Section 4.

\section{The homeomorphism types of the orbit spaces}

In this section we determine the homeomorphism types of those $6$--manifolds
which can appear as the orbit spaces of some regular circle actions on $2$%
--connected $7$--manifolds.

A $2$--connected $7$--manifold $M$ with a regular circle action defines a
principal $S^{1}$--bundle $M\rightarrow N$ with base space $N=M/S^{1}$ \cite%
{Gl}. Fixing an orientation on $S^{1}$ once and for all, and let $N$ be
furnished with the orientation compactible with that on $M$. From the
homotopy exact sequence

\begin{center}
$0\rightarrow \pi _{2}(M)\rightarrow $ $\pi _{2}(N)\rightarrow \pi
_{1}(S^{1})\rightarrow \pi _{1}(M)\rightarrow \pi _{1}(N)\rightarrow 0$
\end{center}

\noindent of the fibration one finds that

\begin{center}
$\pi _{1}(N)=0$; $\pi _{2}(N)\cong \pi _{1}(S^{1})\cong \mathbb{%
\mathbb{Z}
}.$
\end{center}

\noindent Consequently $N$ is a $1$--connected $6$--manifold with $%
H_{2}(N)\cong \mathbb{%
\mathbb{Z}
}$.

Conversely, for a $1$--connected $6$--manifold $N$ with $H_{2}(N)\cong 
\mathbb{%
\mathbb{Z}
}$ let $t\in H^{2}(N)\cong \mathbb{%
\mathbb{Z}
}$ be a generator and let

\begin{center}
$S^{1}\hookrightarrow N_{t}\rightarrow N$
\end{center}

\noindent be the oriented circle bundle over $N$ with Euler class $t$. From
the homotopy exact sequence of this fibration we find that $N_{t}$ is $2$%
--connected with the canonical orientation as the total space of the circle
bundle. Summarizing we get

\bigskip

\noindent \textbf{Lemma 2.1.} Let $S^{1}\times M\rightarrow M$ be a regular
circle action on a $2$--connected $7$--manifold $M$ with orbit space $N$.
Then $N$ is a $1$--connected $6$--manifold with $H_{2}(N)\cong \mathbb{%
\mathbb{Z}
}$.

Conversely, every $1$--connected $6$--manifold $N$ with $H_{2}(N)\cong 
\mathbb{%
\mathbb{Z}
}$ can be realized as the orbit space of some regular circle action on a $2$%
--connected $7$--manifold.$\square $

\bigskip

In view of Lemma 2.1 the classification of those $1$--connected $6$%
--manifolds $N$ with $H_{2}(N)\cong \mathbb{%
\mathbb{Z}
}$ amounts to a crucial step toward a solution to Problem 1.2. In terms of
the known invariants for $1$--connected $6$--manifolds due to Jupp\cite{J}
and Wall\cite{W} we can enumerate all these manifolds in the next result.

Denote by $\Theta $ the set of equivalence classes $[N,t]$ of the pairs $%
\left( N,t\right) $ with $N$ a $1$--connected $6$--manifold whose integral
cohomology satisfies

\begin{center}
$H^{r}(N)=\left\{ 
\begin{tabular}{l}
$\mathbb{%
\mathbb{Z}
}$ if $r=0,2,4,6$; \\ 
$0$ otherwise,%
\end{tabular}%
\right. $
\end{center}

\noindent and with $t\in H^{2}(N)$ a fixed generator. Two elements $%
(N_{1},t_{1}),(N_{2},t_{2})$ are called \textsl{equivalent} if there is a
homeomorphism $f:N_{1}\rightarrow N_{2}$ such that $f^{\ast }t_{2}=t_{1}$.
For each $(N,t)$ fix a generator $x\in H^{4}(N)$ such that the value $%
\left\langle t\cup x,[N]\right\rangle $ of the product $t\cup x$ on the
fundamental class $[N]$ is equal to $1$. Consider the functions

\begin{center}
$k,p:\Theta \rightarrow \mathbb{%
\mathbb{Z}
}$; \ $\varepsilon :\Theta \rightarrow \{0,1\}$; \ $\delta :\Theta
\rightarrow \{0,1\}$
\end{center}

\noindent determined by the following properties

\begin{quote}
i) $t^{2}=k([N,t])x$;

ii) the second Stiefel-Whitney class $w_{2}(N)$ and the first Pontrjagin
class $p_{1}(N)$ of $N$ are given by $\varepsilon ([N,t])t$ $\func{mod}2$
and $p([N,t])x$, respectively;

iii) the class $\Delta (N)\equiv \delta ([N,t])x\func{mod}2\in H^{4}(N;%
\mathbb{%
\mathbb{Z}
}_{2})$ is the Kirby-Siebenmann invariant of $N$,
\end{quote}

\noindent where the Kirby-Siebenmann invariant $\Delta (V)$ of a manifold $V$
is the obstruction to lift the classifying map $V\rightarrow BTOP$ for the
stable tangent bundle of $V$ to $BPL$, and where $BTOP$ and $BPL$ are the
classifying spaces for the stable $TOP$ bundles and $PL$ bundles,
respectively\cite{KS}.

\bigskip

\noindent \textbf{Lemma 2.2.} For each $1$--connected $6$--manifold $M$ with 
$H_{2}(M)\cong \mathbb{%
\mathbb{Z}
}$ there exists an $r\in \mathbb{%
\mathbb{N}
}$ and an element $[N,t]\in \Theta $ such that $M\cong \#_{r}S^{3}\times
S^{3}\#N$.

Moreover, the system $\{k,p,\varepsilon ,\delta \}$ is a set of complete
invariants for elements $\left[ N,t\right] \in \Theta $ that is subject to
the following constraints: \ \ \ \ \ \ \ \ \ \ \ \ \ \ \ \ 

\begin{quote}
i)If $k(\left[ N,t\right] )\equiv 1\func{mod}2$, then $\varepsilon (\left[
N,t\right] )=0$ and

$p(\left[ N,t\right] )=24m+4k(\left[ N,t\right] )+24\delta (\left[ N,t\right]
)$ for some $m\in 
\mathbb{Z}
$;

ii)If $k(\left[ N,t\right] )\equiv 0\func{mod}2$, then for some $m\in 
\mathbb{Z}
$

$p(\left[ N,t\right] )=\left\{ 
\begin{array}{cc}
24m+4k(\left[ N,t\right] )+24\delta (\left[ N,t\right] ) & \text{if }%
\varepsilon (\left[ N,t\right] )=0 \\ 
48m+k(\left[ N,t\right] )+24\delta (\left[ N,t\right] ) & \text{if }%
\varepsilon (\left[ N,t\right] )=1%
\end{array}%
\right. $.
\end{quote}

\noindent In addition, the manifold $N$ is smoothable if and only if $\delta
(\left[ N,t\right] )=0.$

\bigskip

\noindent \textbf{Proof.} This is a direct consequence of \cite[Theorem 0;
Theorem 1]{J}. In particular, the expressions of the function $p$ are
deduced from the following relation on $H^{6}(N)$ which holds for all $d\in 
\mathbb{%
\mathbb{Z}
}$:

\begin{center}
$(2dt+\varepsilon (\left[ N,t\right] )t)^{3}\equiv (p(\left[ N,t\right]
)x+24\delta (\left[ N,t\right] )x)(2dt+\varepsilon (\left[ N,t\right] )t)%
\func{mod}48$.$\square $
\end{center}

\section{Circle bundles over $\left[ N,t\right] \in \Theta $}

Lemma 2.2 singles out the family $\Theta $ of $1$--connected $6$--manifolds
which plays a key role in presenting the orbit spaces of regular circle
actions on $2$--connected $7$--manifolds. In this section we determine the
homeomorphism and diffeomorphism type of the total space $N_{t}$ of the
circle bundle over $\left[ N,t\right] \in \Theta $, i.e. the oriented circle
bundle over $N$ with Euler class $t$. For this purpose we shall recall in
Section 3.1 that the definition of the known invariant system for $2$%
--connected $7$--manifolds. The main results in this section are Lemmas 3.3
and 3.5, which identify the homeomorphism and diffeomorphism types of the
manifolds $N_{t}$ with certain $M_{l,k}^{c}$.

\subsection{Invariants for 2--connected 7--manifolds}

Recall from Eells, Kuiper\cite{E}, Kreck, Stolz\cite{K} and Wilkens\cite{Wl}
that associated to each $2$--connected $7$--manifold $M$ there is a system $%
\{H,\frac{p_{1}}{2},b,\Delta ,\mu ,s_{1}\}$ of invariants characterized by
the following properties:

\begin{quote}
i) $H$ is the forth integral cohomology group $H^{4}(M)$\cite{Wl};

ii) $\frac{p_{1}}{2}(M)\in H$ is the first spin characteristic class \cite%
{Th} introduced by Wilkens \cite{Wl} in smooth category and extended to
topological category by Kreck--Stolz\cite[Lemma 6.5]{K};

iii) $b:\tau (H)\otimes \tau (H)\rightarrow \mathbb{%
\mathbb{Q}
}/\mathbb{%
\mathbb{Z}
}$ is the linking form on the torsion part $\tau (H)$ of the group $H$\cite%
{Wl};

iv) $\Delta (M)\in H^{4}(M;\mathbb{%
\mathbb{Z}
}_{2})$ is the Kirby--Siebenmann invariant of $M$ \cite{K}.
\end{quote}

\noindent Furthermore, if the manifold $M$ is smooth and bounds a smooth $8$%
--manifold $W$ with the induced map $j^{\ast }:H^{4}(W,M;\mathbb{%
\mathbb{Q}
})\rightarrow H^{4}(W;\mathbb{%
\mathbb{Q}
})$ an isomorphism, then

\begin{quote}
v) the invariant $\mu \in \mathbb{%
\mathbb{Q}
}/\mathbb{%
\mathbb{Z}
}$ is firstly defined in \cite{E} for a spin $W$ and extended in \cite{K}
for a general $W$, whose value is given by the formula

$\mu (M)\equiv $--$\frac{1}{2^{5}\cdot 7}\sigma (W)+\frac{1}{2^{7}\cdot 7}%
p_{1}^{2}(W)$--$\frac{1}{2^{6}\cdot 3}z^{2}\cdot p_{1}(W)+\frac{1}{%
2^{7}\cdot 3}z^{4}\func{mod}%
\mathbb{Z}
$,
\end{quote}

\noindent where $z\in H^{2}(W)$ satisfies $w_{2}(W)=z\func{mod}2$, $\sigma
(W)$ is the signature of the intersection form on $H^{4}(W,M;%
\mathbb{Q}
)$, and where $p_{1}^{2}(W)$, $z^{2}\cdot p_{1}(W)$ and $z^{4}$ are the
characteristic numbers

\begin{quote}
$\langle p_{1}(W)\cup j^{\ast -1}p_{1}(W),[W,M]\rangle $, $\langle z^{2}\cup
j^{\ast -1}p_{1}(W),[W,M]\rangle $,

$\langle z^{2}\cup j^{\ast -1}z^{2},[W,M]\rangle $,
\end{quote}

\noindent respectively. Finally, if $M$ is topological and bounds a
topological 8--manifold $W$ with the induced map $j^{\ast }:H^{4}(W,M;%
\mathbb{%
\mathbb{Q}
})\rightarrow H^{4}(W;\mathbb{%
\mathbb{Q}
})$ an isomorphism, then

\begin{quote}
vi) the topological invariant $s_{1}\in \mathbb{%
\mathbb{Q}
}/\mathbb{%
\mathbb{Z}
}$ is defined in \cite{K} whose value is given by

$s_{1}(M)\equiv $--$\frac{1}{2^{3}}\sigma (W)+\frac{1}{2^{5}}p_{1}^{2}(W)$--$%
\frac{7}{2^{4}\cdot 3}z^{2}\cdot p_{1}(W)+\frac{7}{2^{5}\cdot 3}z^{4}\func{%
mod}%
\mathbb{Z}
$.
\end{quote}

\noindent \textbf{Example 3.1.} Let $N_{t}$ be the total space of the circle
bundle over $\left[ N,t\right] \in \Theta $. Then the system $\{H,\frac{p_{1}%
}{2},b,\Delta ,\mu ,s_{1}\}$ of invariants for the manifold $N_{t}$ can be
expressed in terms of the invariants for $\Theta $ introduced in Lemma 2.2
as follows. For simplicity we write $p,k,\varepsilon $ and $\delta $ in
place of $p(\left[ N,t\right] )$, $k(\left[ N,t\right] )$, $\varepsilon (%
\left[ N,t\right] )$ and $\delta (\left[ N,t\right] )$, respectively.

\begin{quote}
i) $H^{4}(N_{t})\cong \mathbb{%
\mathbb{Z}
}_{k}$ with generator $\pi ^{\ast }(x)$, where $\pi :N_{t}\rightarrow N$ is
the bundle projection and

$\mathbb{%
\mathbb{Z}
}_{k}=\left\{ 
\begin{array}{cc}
\mathbb{Z}
& \text{if }k=0 \\ 
\mathbb{Z}
/k%
\mathbb{Z}
& \text{if }k\neq 0%
\end{array}%
\right. $;

ii) $\Delta (N_{t})\equiv \frac{1+(-1)^{k}}{2}\cdot \delta \pi ^{\ast }(x)%
\func{mod}2$;

iii) $\frac{p_{1}}{2}(N_{t})\equiv \frac{p+\varepsilon k}{2}\pi ^{\ast }(x)%
\func{mod}k$;

iv) $b(\pi ^{\ast }(x),\pi ^{\ast }(x))\equiv \frac{1}{k}\func{mod}\mathbb{%
\mathbb{Z}
}$;

v) $\mu (N_{t})\equiv $--$\frac{\left\vert k\right\vert }{2^{5}\cdot 7k}+%
\frac{(p+k)^{2}}{2^{7}\cdot 7k}+\frac{(\varepsilon -1)(2p+k)}{2^{7}\cdot 3}%
\func{mod}%
\mathbb{Z}
$;

vi) $s_{1}(N_{t})=$--$\frac{\left\vert k\right\vert }{2^{3}k}+\frac{(p+k)^{2}%
}{2^{5}k}+\frac{7(\varepsilon -1)(2p+k)}{2^{5}\cdot 3}\func{mod}%
\mathbb{Z}
$.
\end{quote}

Firstly, from the section $H^{2}(N)\overset{\cup t}{\rightarrow }H^{4}(N)%
\overset{\pi ^{\ast }}{\rightarrow }H^{4}(N_{t})\rightarrow 0$ in the Gysin
sequence of the fibration $N_{t}\overset{\pi }{\rightarrow }N$ and from the
relation $t^{2}=kx$ on $H^{4}(N)$ we find that $H^{4}(N_{t})\cong \mathbb{%
\mathbb{Z}
}_{k}$ with generator $\pi ^{\ast }(x)$. This shows i).

Next, let $f:N\rightarrow BTOP$ be the classifying map for the stable
tangent bundle of $N$. In view of the decomposition $TN_{t}\cong \pi ^{\ast
}TN\oplus \varepsilon ^{1}$($\varepsilon ^{1}$ denotes the trivial line
bundle) for the tangent bundle of $N_{t}$ the classifying map for the stable
tangent bundle of $N_{t}$ is given by the composition $f\circ \pi
:N_{t}\rightarrow N\rightarrow BTOP$. It follows that the Kirby-Siebenmann
invariant $\Delta (N_{t})$ of the manifold $N_{t}$ is $\pi ^{\ast }\Delta
(N)\equiv \delta \pi ^{\ast }(x)\func{mod}2$. This shows ii).

To calculate the remaining invariants $\frac{p_{1}}{2},b,\mu ,s_{1}$ of the
manifold $N_{t}$ we make use of the associated disk bundle $W_{t}\overset{%
\pi _{0}}{\rightarrow }N$ of the oriented $2$--plane bundle $\xi _{t}$ over $%
N$ with Euler class $t$. If $\varepsilon =1$, It follows from the
decomposition $TW_{t}\cong \pi _{0}^{\ast }TN\oplus \pi _{0}^{\ast }\xi _{t}$
that $w_{2}(W_{t})=0$ and $\frac{p_{1}}{2}(W_{t})=\frac{p+k}{2}\pi
_{0}^{\ast }x.$ From the relation $\partial W_{t}=N_{t}$ we get

\begin{quote}
$\frac{p_{1}}{2}(N_{t})=\frac{p+k}{2}\pi ^{\ast }x$.
\end{quote}

\noindent If $\varepsilon =0$, from the decomposition $TN_{t}\cong \pi
^{\ast }TN\oplus \varepsilon ^{1}$ we get

\begin{quote}
$\frac{p_{1}}{2}(N_{t})=\frac{p}{2}\pi ^{\ast }(x)$.
\end{quote}

\noindent This shows iii).

To compute the linking form $b$ of $N_{t}$ we can assume that $k\neq 0$.
Consider the commutative ladder of exact sequences

\begin{center}
$%
\begin{array}{ccccccccc}
0 & \rightarrow & H^{4}(W_{t},N_{t}) & \overset{j^{\ast }}{\rightarrow } & 
H^{4}(W_{t}) & \overset{i^{\ast }}{\rightarrow } & H^{4}(N_{t}) & \rightarrow
& 0 \\ 
&  & \cong \uparrow \phi &  & \uparrow \pi _{0}^{\ast } &  & \uparrow id & 
&  \\ 
0 & \rightarrow & H^{2}(N) & \overset{\cup t}{\rightarrow } & H^{4}(N) & 
\overset{\pi ^{\ast }}{\rightarrow } & H^{4}(N_{t}) & \rightarrow & 0%
\end{array}%
$
\end{center}

\noindent with $\phi $ the Thom isomorphism. Since $\pi ^{\ast }(x)=i^{\ast
}\pi _{0}^{\ast }(x)$ and $y:=\phi (t)$ is a generator of $%
H^{4}(W_{t},N_{t}) $ with $j^{\ast }(y)=\pi _{0}^{\ast }(t^{2})=k\pi
_{0}^{\ast }(x)$ we get

\begin{quote}
$b(\pi ^{\ast }(x),\pi ^{\ast }(x))\equiv \frac{1}{k}<y\cup \pi _{0}^{\ast
}x,[W_{t},N_{t}]>\equiv \frac{1}{k}\func{mod}\mathbb{%
\mathbb{Z}
}$.
\end{quote}

\noindent This shows iv).

Since the induced map $j^{\ast }:H^{4}(W_{t},N_{t};\mathbb{%
\mathbb{Q}
})\rightarrow H^{4}(W_{t};\mathbb{%
\mathbb{Q}
})$ is clearly an isomorphism when $k\neq 0$, the invariants $\mu $ and $%
s_{1}$ are defined for $N_{t}$. Moreover, from the Lefschetz duality and the
relation $j^{\ast }(y)=k\pi _{0}^{\ast }(x)$ we get $\sigma (W_{t})=\frac{%
\left\vert k\right\vert }{k}.$ From the decomposition $TW_{t}\cong \pi
_{0}^{\ast }(TN\oplus \xi _{t})$ we get, in addition to

\begin{quote}
$p_{1}(W_{t})=\pi _{0}^{\ast }(p_{1}(N)+t^{2})=(p+k)\pi _{0}^{\ast }(x)$,
\end{quote}

\noindent that

\begin{quote}
$w_{2}(W_{t})\equiv (\varepsilon +1)\pi _{0}^{\ast }(t)\func{mod}2$.
\end{quote}

\noindent Therefore we can take $z=(1$--$\varepsilon )\pi _{0}^{\ast }(t)$
in the formulae for $\mu $ and $s_{1}$, and as a result

\begin{quote}
$z^{2}=(1$--$\varepsilon )^{2}\pi _{0}^{\ast }(t^{2})=k(1$--$\varepsilon
)^{2}\pi _{0}^{\ast }(x)$.
\end{quote}

\noindent As the group $H^{4}(W_{t},N_{t})\cong \mathbb{%
\mathbb{Z}
}$ is generated by $y=\phi (t)$ with the relation $j^{\ast }(y)=k\pi
_{0}^{\ast }(x),$ the isomorphism

\begin{quote}
$H^{4}(W_{t},N_{t})\otimes H^{4}(W_{t})\overset{\cup }{\rightarrow }%
H^{8}(W_{t},N_{t})$
\end{quote}

\noindent by the Lefschetz duality, together with the formulae for $%
p_{1}(W_{t})$ and $z^{2}$ above, implies the relations below

\begin{quote}
$z^{2}p_{1}(W_{t})=(1$--$\varepsilon )^{2}(p+k)$; $p_{1}^{2}(W_{t})=\frac{1}{%
k}(p+k)^{2}$; $z^{4}=k(1$--$\varepsilon )^{4}$.
\end{quote}

\noindent Substituting these values in the formulae for $\mu $ and $s_{1}$
yields v) and vi) respectively. This completes the computation of the
invariant system for the manifolds $N_{t}$.

\bigskip

\noindent \textbf{Example} \textbf{3.2}. The invariant system $\{H,\Delta ,%
\frac{p_{1}}{2},b,s_{1},\mu \}$ of the manifolds $M_{l,k}^{c}$ has been
computed by Crowley and Escher \cite{C} for the case of $c=0$. We extend
their calculation as to include the exceptional case of $c=1$.

\begin{quote}
i) $H^{4}(M_{l,k}^{c})\cong \mathbb{Z}_{k}$ with generator $\kappa =\left\{ 
\begin{array}{cc}
\pi _{M}^{\ast }(\iota ) & \text{if }c=0 \\ 
(\pi _{M}\circ h_{M})^{\ast }(\iota ) & \text{if }c=1%
\end{array}%
\right. $;

ii) $b(\kappa ,\kappa )\equiv \frac{1}{k}\func{mod}\mathbb{Z}$;

iii) $\Delta (M_{l,k}^{c})\equiv \frac{1+(-1)^{k}}{2}\cdot c\kappa \func{mod}%
2$;

iv) $\frac{p_{1}}{2}(M_{l,k}^{c})\equiv (2l+12c)\kappa \func{mod}k$;

v) $s_{1}(M_{l,k}^{c})\equiv \frac{(2l+k+12c)^{2}-\left\vert k\right\vert }{%
8k}\func{mod}\mathbb{Z}$;

vi) $\mu (M_{l,k}^{0})\equiv \frac{(k+2l)^{2}-\left\vert k\right\vert }{%
28\cdot 8k}\func{mod}\mathbb{Z}$.
\end{quote}

Firstly, since $h_{M}:M_{l,k}^{1}\rightarrow $ $M_{l,k}^{0}$ is a homotopy
equivalence we get i) and ii) from the relations $H^{4}(M_{l,k}^{0})\cong 
\mathbb{Z}_{k}$ (with generator $\pi _{M}^{\ast }(\iota )$) and $b(\pi
_{M}^{\ast }(\iota ),\pi _{M}^{\ast }(\iota ))\equiv \frac{1}{k}\func{mod}%
\mathbb{Z}$ when $c=0$.

Next, since the map $\mathcal{S}^{TOP}(M_{l,k}^{0})\overset{\Delta }{%
\rightarrow }H^{4}(M_{l,k}^{0};\mathbb{Z}_{2})$ of taking Kirby-- Siebenmann
class is a surjective homomorphism \cite[Theorem 15.1]{S}, and since $\left[
M_{l,k}^{1},h_{M}\right] $ is a generator of the cyclic group $\mathcal{S}%
^{TOP}(M_{l,k}^{0})\cong \mathbb{Z}_{k}$ we have

\begin{quote}
$\Delta (M_{l,k}^{1})=\Delta (\left[ M_{l,k}^{1},h_{M}\right] )=\frac{%
1+(-1)^{k}}{2}\kappa \func{mod}2$.
\end{quote}

\noindent This shows iii).

To calculate the remaining invariants $\frac{p_{1}}{2},\mu ,s_{1}$ of the
manifold $M_{l,k}^{c}$ we construct an 8--manifold $W_{l,k}^{c}$ with
boundary $\partial W_{l,k}^{c}\cong M_{l,k}^{c}$ as follows. Let $\pi
_{W}:W_{l,k}^{0}\rightarrow S^{4}$ be the associated disk bundle of the
sphere bundle $\pi _{M}:M_{l,k}^{0}\rightarrow S^{4}$. Then $%
M_{l,k}^{0}=\partial W_{l,k}^{0}$. Write $\mathcal{S}^{TOP}(W_{l,k}^{0})$
for the set of equivalence classes $[W,h]$ of the pairs $(W,h)$ with $%
h:(W,\partial W)\rightarrow (W_{l,k}^{0},M_{l,k}^{0})$ a homotopy
equivalence between $8$-manifolds with boundary. Consider the following
commutative diagram analoguing to the one \cite[(7)]{C} due to Crowley and
Escher in the PL--category(see also Section 1)

\begin{center}
$%
\begin{array}{ccccc}
\mathcal{S}^{TOP}(W_{l,k}^{0}) & \overset{\eta }{\underset{\cong }{%
\rightarrow }} & [W_{l,k}^{0},G/TOP] & \underset{\cong }{\overset{d}{%
\rightarrow }} & H^{4}(W_{l,k}^{0})\cong 
\mathbb{Z}
\\ 
i^{\ast }\downarrow &  & i^{\ast }\downarrow &  & i^{\ast }\downarrow \\ 
\mathcal{S}^{TOP}(M_{l,k}^{0}) & \overset{\eta }{\underset{\cong }{%
\rightarrow }} & [M_{l,k}^{0},G/TOP] & \underset{\cong }{\overset{d}{%
\rightarrow }} & H^{4}(M_{l,k}^{0})\cong 
\mathbb{Z}
_{k}%
\end{array}%
$
\end{center}

\noindent where $i^{\ast }:\mathcal{S}^{TOP}(W_{l,k}^{0})\rightarrow 
\mathcal{S}^{TOP}(M_{l,k}^{0})$ sends each $\left[ W,h\right] $ to the
restriction $\left[ \partial W,h|_{\partial W}\right] $. Writing $\left[
W_{l,k}^{1},h_{W}\right] $ for the element $(d\circ \eta )^{-1}(\pi
_{W}^{\ast }(\iota ))\in \mathcal{S}^{TOP}(W_{l,k}^{0})$ we get $%
M_{l,k}^{1}\cong \partial W_{l,k}^{1}$ from the diagram above.

To find the formula of $\frac{p_{1}}{2}(M_{l,k}^{c})$ we compute the first
Pontrjagin class $p_{1}(W_{l,k}^{c})$ of $W_{l,k}^{c}$. Let $\alpha $ denote
the generator of $H^{4}(W_{l,k}^{c})\cong 
\mathbb{Z}
$ satisfies

\begin{center}
$\alpha =\left\{ 
\begin{array}{cc}
\pi _{W}^{\ast }(\iota ) & \text{if }c=0\text{;} \\ 
(\pi _{W}\circ h_{W})^{\ast }(\iota ) & \text{if }c=1\text{;}%
\end{array}%
\right. $
\end{center}

\noindent and associate an integer $p(W_{l,k}^{c})$ to $W_{l,k}^{c}$ such
that $p_{1}(W_{l,k}^{c})=p(W_{l,k}^{c})\alpha $. Let $\overline{i}%
:G/TOP\rightarrow BTOP$ be the natural inclusion and let $%
f_{c}:W_{l,k}^{c}\rightarrow BTOP$ be the classifying map for the stable
tangent bundle of $W_{l,k}^{c}$. It follows from the isomorphism $\mathcal{S}%
^{TOP}(W_{l,k}^{0})\overset{\eta }{\rightarrow }[W_{l,k}^{0},G/TOP]$ and the
proof of \cite[Theorem 2.23]{Ma} that

\begin{quote}
$\overline{i}_{\ast }\eta ([W_{l,k}^{1},h_{W}])=h_{W}^{\ast -1}[f_{1}]$--$%
[f_{0}]$
\end{quote}

\noindent and hence

\begin{quote}
$p(W_{l,k}^{1})\pi _{W}^{\ast }(\iota )=h_{W}^{\ast
-1}p_{1}(W_{l,k}^{1})=p_{1}(W_{l,k}^{0})+f^{\ast }\overline{i}^{\ast }p_{1}$,
\end{quote}

\noindent where $f$ $=d^{-1}(\pi _{W}^{\ast }(\iota ))=\eta
([W_{l,k}^{1},h_{W}])$ is the generator of $[W_{l,k}^{0},G/TOP]$ and $%
p_{1}\in H^{4}(BTOP)$ is the first Pontrjagin class \cite{J}. It is shown in 
\cite[Lemma 13.3,Proposition 13.4]{S} that a generator $g$ of $[S^{4},G/TOP]$
corresponds to a topological bundle $\xi $ with classifying map $\overline{i}%
\circ g$ and Pontrjagin class $p_{1}(\xi )=g^{\ast }\overline{i}^{\ast
}p_{1}=\pm 24\iota $. With an appropriate choice of $\pm
d:[W_{l,k}^{0},G/TOP]\rightarrow H^{4}(W_{l,k}^{0})$ applying $\pi
_{W}^{\ast }$ to this equation we get

\begin{quote}
$f^{\ast }\overline{i}^{\ast }p_{1}=24\pi _{W}^{\ast }(\iota ).$
\end{quote}

\noindent This, together with the fact $p(W_{l,k}^{0})=2(k+2l)$ \cite{M} and
the formula for $p(W_{l,k}^{1})$ above, implies that $%
p(W_{l,k}^{c})=2k+4l+24c$. Consequently from $M_{l,k}^{c}\cong \partial
W_{l,k}^{c}$ we get iv).

Finally, we compute $s_{1}(M_{l,k}^{c})$. The exact sequence

\begin{quote}
$H^{4}(W_{l,k}^{c},M_{l,k}^{c})\overset{j^{\ast }}{\rightarrow }%
H^{4}(W_{l,k}^{c})\rightarrow H^{4}(M_{l,k}^{c})\rightarrow 0$,
\end{quote}

\noindent together with the isomorphisms $H^{4}(M_{l,k}^{c})\cong \mathbb{Z}%
_{k}$ and $H^{4}(W_{l,k}^{c},M_{l,k}^{c})\cong \mathbb{Z}$ by the Lefschetz
duality, implies that we can take a generator $\beta $ of $%
H^{4}(W_{l,k}^{c},M_{l,k}^{c})$ such that $j^{\ast }(\beta )=k\alpha $.
Since $j^{\ast }:H^{4}(W_{l,k}^{c},M_{l,k}^{c};\mathbb{Q})\rightarrow
H^{4}(W_{l,k}^{c};\mathbb{Q})$ is an isomorphism for $k\neq 0$ the invariant 
$s_{1}$ is defined for $M_{l,k}^{c}$. It follows from the Lefschetz duality
and the relation $j^{\ast }(\beta )=k\alpha $ that $\sigma (W_{l,k}^{c})=%
\frac{\left\vert k\right\vert }{k}$. On the other hand, the formula for $%
p(W_{l,k}^{c})$, together with the relation $j^{\ast }(\beta )=k\alpha $ and
the isomorphism

\begin{quote}
$H^{4}(W_{l,k}^{c},M_{l,k}^{c})\otimes H^{4}(W_{l,k}^{c})\overset{\cup }{%
\rightarrow }H^{8}(W_{l,k}^{c},M_{l,k}^{c})$
\end{quote}

\noindent by the Lefschetz duality, implies that

\begin{quote}
$p_{1}^{2}(W_{l,k}^{c})=\frac{4}{k}(k+2l+12c)^{2}$.
\end{quote}

\noindent In addition, the relation $w_{2}(W_{l,k}^{1})=w_{2}(W_{l,k}^{0})=0$
indicates that we can take $z=0$ in the formula of $s_{1}$. Substituting the
values of $\sigma (W_{l,k}^{c}),z,p_{1}^{2}(W_{l,k}^{c})$ in the formula for 
$s_{1}$ shows v).

Similarly, we refer vi) to Crowley and Escher \cite{C}. This completes the
computation of the invariant system for the manifolds $M_{l,k}^{c}$.

\subsection{Circle bundles over $[N,t]\in \Theta $}

In this section we will prove Lemmas 3.3 and 3.5 which identify the
homeomorphism and diffeomorphism types of the manifolds $N_{t}$ with certain 
$M_{l,k}^{c}$.

\bigskip

\noindent \textbf{Lemma 3.3.} Let $N_{t}$ be the total space of the circle
bundle over $[N,t]\in \Theta $. Then there is a homeomorphism $N_{t}\cong
M_{l,k}^{c}$ where

\begin{quote}
$(k,c)=(k(\left[ N,t\right] ),\frac{1+(-1)^{k(\left[ N,t\right] )}}{2}\cdot
\delta (\left[ N,t\right] ));$

$l=\frac{p(\left[ N,t\right] )+(3\varepsilon (\left[ N,t\right] )-4)\cdot k(%
\left[ N,t\right] )-(1+(-1)^{k(\left[ N,t\right] )})\cdot 12\delta (\left[
N,t\right] )}{4}.$
\end{quote}

\noindent \textbf{Proof.} We divide the proof into two cases depending on $%
\Delta (N_{t})\equiv 0$ or $1$ $\func{mod}2$.

Case 1. $\Delta (N_{t})\equiv 0\func{mod}2$ (i.e. the manifold $N_{t}$ is
smoothable, see \cite[p.33]{Ma}, \cite{HM} and \cite[Theorem 5.4]{S}):
Module by a $%
\mathbb{Z}
_{2}$ ambiguity Wilkens\cite{Wl} showed that the system $\{H,\frac{p_{1}}{2}%
,b\}$ of invariants classifies $N_{t}$ and $M_{l,k}^{0}$ up to
homeomorphism. Moreover, Crowley and Escher \cite{C} proved that this
ambiguity can be realized by some $M_{l,k}^{0}$ whose homeomorphism types
can be distinguished by the invariant $s_{1}$ and hence the system $\{H,%
\frac{p_{1}}{2},b,s_{1}\}$ classifies the manifolds $N_{t}$ and $M_{l,k}^{0}$%
. Therefore the proof is completed by comparing these invariants for $N_{t}$
and $M_{l,k}^{0}$ obtained in Example 3.1 and 3.2.

Case 2. $\Delta (N_{t})\equiv 1\func{mod}2$: We only need to show that $%
N_{t} $ is homeomorphic to $M_{l,k}^{1}$ with $\left[ N,t\right] \in \Theta $
and

\begin{quote}
$(k,l)=(k(\left[ N,t\right] ),\frac{p(\left[ N,t\right] )+(3\varepsilon (%
\left[ N,t\right] )-4)\cdot k(\left[ N,t\right] )-24}{4})$.
\end{quote}

\noindent It suffices to construct a homotopy equivalence $%
q:N_{t}\rightarrow M_{l,k}^{0}$ with

\begin{quote}
$\left[ N_{t},q\right] =\left[ M_{l,k}^{1},h_{M}\right] \in \mathcal{S}%
^{TOP}(M_{l,k}^{0})$.
\end{quote}

According to Lemma 2.2 there exists a manifold $N^{\prime }$ with $%
[N^{\prime },t^{\prime }]\in \Theta $ whose invariant system $(k([N^{\prime
},t^{\prime }]),p([N^{\prime },t^{\prime }]),\varepsilon ([N^{\prime
},t^{\prime }]),\delta ([N^{\prime },t^{\prime }])$ is

\begin{quote}
$(k([N,t]),p(\left[ N,t\right] )$--$24,\varepsilon ([N,t]),0)$.
\end{quote}

\noindent Consider the map $\eta :\mathcal{S}^{TOP}(N^{^{\prime
}})\rightarrow \lbrack N^{^{\prime }},G/TOP]$ in the surgery exact sequence
of $N^{^{\prime }}$. By the argument at the end of the proof of \cite[%
Theorem 1]{J} we find a homotopy equivalence $h_{N}:N\rightarrow N^{^{\prime
}}$ such that

\begin{quote}
i) the class $\eta (\left[ N,h_{N}\right] )$ is trivial on the $2$ skeleton
of $N^{^{\prime }}$;

ii) the primary obstruction to finding a null-homotopy of $\eta (\left[
N,h_{N}\right] )$ is the generator $x^{\prime }\in H^{4}(N^{^{\prime }};\pi
_{4}(G/TOP))$ with $\left\langle x^{\prime }\cup t^{\prime },[N^{\prime
}]\right\rangle $=$1$.
\end{quote}

\noindent Pulling $h_{N}$ back by the bundle projection $\pi ^{\prime
}:N_{t}^{^{\prime }}\rightarrow N^{^{\prime }}$ induces a homotopy
equivalence $h_{t}:N_{t}\rightarrow N_{t}^{^{\prime }}$. On the other hand,
by the result of Case 1 we get a homeomorphism $u:M_{l,k}^{0}\rightarrow
N_{t}^{^{\prime }}$ such that $u^{\ast }(\pi ^{\prime \ast }(x^{\prime
}))=\pi _{M}^{\ast }(\iota )$. So it remains to show that $\left[
N_{t},u^{-1}\circ h_{t}\right] =\left[ M_{l,k}^{1},h_{M}\right] $.

Let $[N^{^{\prime }},G/TOP]_{2}$ denote the subset of $[N^{^{\prime
}},G/TOP] $ whose elements are trivial on the $2$ skeleton of $N^{^{\prime
}} $ and consider the following two commutative diagrams:

\begin{center}
$%
\begin{array}{ccc}
\mathcal{S}^{TOP}(N^{^{\prime }}) & \overset{\pi ^{\prime \ast }}{%
\rightarrow } & \mathcal{S}^{TOP}(N_{t}^{^{\prime }}) \\ 
\downarrow \eta &  & \cong \downarrow \eta \\ 
\left[ N^{^{\prime }},G/TOP\right] & \overset{\pi ^{\prime \ast }}{%
\rightarrow } & \left[ N_{t}^{^{\prime }},G/TOP\right]%
\end{array}%
\quad $

\bigskip

$%
\begin{array}{ccccc}
&  & \mathcal{S}^{TOP}(N_{t}^{^{\prime }}) & \underset{\cong }{\overset{%
u^{\ast }}{\rightarrow }} & \mathcal{S}^{TOP}(M_{l,k}^{0}) \\ 
&  & \cong \downarrow \eta &  & \cong \downarrow \eta \\ 
\left[ N^{^{\prime }},G/TOP\right] _{2} & \overset{\pi ^{\prime \ast }}{%
\rightarrow } & \left[ N_{t}^{^{\prime }},G/TOP\right] & \underset{\cong }{%
\overset{u^{\ast }}{\rightarrow }} & \left[ M_{l,k}^{0},G/TOP\right] \\ 
\downarrow d &  & \cong \downarrow d &  & \cong \downarrow d \\ 
H^{4}(N^{^{\prime }}) & \overset{\pi ^{\prime \ast }}{\rightarrow } & 
H^{4}(N_{t}^{^{\prime }}) & \underset{\cong }{\overset{u^{\ast }}{%
\rightarrow }} & H^{4}(M_{l,k}^{0})%
\end{array}%
$
\end{center}

\noindent where

\begin{quote}
i) $\mathcal{S}^{TOP}(N^{^{\prime }})\overset{\pi ^{\prime \ast }}{%
\rightarrow }\mathcal{S}^{TOP}(N_{t}^{^{\prime }})$ maps $\left[ N^{^{\prime
\prime }},h^{^{\prime \prime }}\right] $ to $\left[ N_{t}^{^{\prime \prime
}},h_{t}^{^{\prime \prime }}\right] $ with $h_{t}^{^{\prime \prime }}$ a
pull-back of $h^{^{\prime \prime }}$ by the bundle projection $\pi ^{\prime
}:N_{t}^{^{\prime }}\rightarrow N^{^{\prime }};$

ii) $\mathcal{S}^{TOP}(N_{t}^{^{\prime }})\overset{u^{\ast }}{\rightarrow }$ 
$\mathcal{S}^{TOP}(M_{l,k}^{0})$ maps $\left[ M^{^{\prime }},g^{^{\prime }}%
\right] $ to $\left[ M^{^{\prime }},u^{-1}\circ g^{^{\prime }}\right] ;$

iii) the maps $d$ send a homotopy class to its primary obstruction to
null-homotopy.
\end{quote}

\noindent The diagrams above, together with the relations

\begin{quote}
$\pi ^{\prime \ast }\left[ N,h_{N}\right] =\left[ N_{t},h_{t}\right] $, $%
u^{\ast }(\pi ^{\prime \ast }(x^{\prime }))=\pi _{M}^{\ast }(\iota )$ and $%
d(\eta \left[ N,h_{N}\right] )=x^{\prime }$,
\end{quote}

\noindent imply that $u^{\ast }\left[ N_{t},h_{t}\right] =\left[
M_{l,k}^{1},h_{M}\right] $, i.e. $\left[ N_{t},u^{-1}\circ h_{t}\right] =%
\left[ M_{l,k}^{1},h_{M}\right] $. This completes the proof of Case 2.$%
\square $

\bigskip

The following lemma plays a key role in the proof of Lemma 3.5\ and its
proof will be postponed to the end of this section.

\bigskip

\noindent \textbf{Lemma 3.4}. Let $N_{t}$ be the total space of the circle
bundle over $[N,t]\in \Theta $ with $\delta ([N,t])=0$ and $k([N,t])=0$.
Then there exists an $8$--manifold $W$ homotopy equivalent to $S^{4}$ whose
boundary satisfies

\begin{quote}
$\partial W\cong \left\{ 
\begin{array}{cc}
N_{t} & \text{if }\varepsilon ([N,t])=1 \\ 
N_{t}\#\Sigma _{\frac{p([N,t])}{24}} & \text{if }\varepsilon ([N,t])=0%
\end{array}%
\right. .$
\end{quote}

\noindent \textbf{Lemma 3.5}. Let $N_{t}$ be the total space of the circle
bundle over $[N,t]\in \Theta $ with $\delta ([N,t])=0$. Then one has a
diffeomorphism $N_{t}\cong M_{l,k}^{0}\#\Sigma _{r}$ where $N_{t}$ has the
smooth structure as the total space of the circle bundle and where

\begin{center}
$(k,l,r)=(k(\left[ N,t\right] ),\frac{p(\left[ N,t\right] )+(3\varepsilon (%
\left[ N,t\right] )-4)\cdot k(\left[ N,t\right] )}{4},\frac{(1-\varepsilon (%
\left[ N,t\right] ))\cdot (p(\left[ N,t\right] )-4k(\left[ N,t\right] ))}{24}%
)$.
\end{center}

\noindent \textbf{Proof.} In the case of $k([N,t])\neq 0$ it is shown in 
\cite{C} that the system $\{H,\frac{p_{1}}{2},b,\mu \}$ classifies $N_{t}$
and $M_{l,k}^{0}$ up to diffeomorphism. Hence the proof is done by comparing
those invariants for $N_{t}$ and $M_{l,k}^{0}$ obtained in Examples 3.1 and
3.2.

Assume next that $k([N,t])=0$ and let $W$ be the $8$--manifold given by
Lemma 3.4. We can take a closed tubular neighborhood $E$ of an embedding $%
h:S^{4}\hookrightarrow $ Interior $W$ which is also a homotopy equivalence 
\cite[Lemma 6]{MM}. As $H_{i}(W\backslash $Interior $E,\partial E)\cong
H_{i}(W,E)=0$ for all $i$ by the excision theorem, $W\backslash $Interior $E$
is an h--cobordism between $\partial W$ and $\partial E$. Hence we get the
diffeomorphisms $\partial W\cong \partial E\cong M_{l,k}^{0}$ for some $%
l,k\in 
\mathbb{Z}
$ by the h-cobordism theorem\cite[Theorem 9.1]{M3} and the fact that $E$ is
the total space of the normal disk bundle of the embedding $h$. Comparing
the invariants $\{H,\frac{p_{1}}{2}\}$ of $\partial W$ and $M_{l,k}^{0}$
given in Examples 3.1 and 3.2 we find that $\partial W$ is diffeomorphic to $%
M_{\frac{p([N,t])}{4},0}^{0}$\cite{C}. This shows

\begin{quote}
$N_{t}\cong \left\{ 
\begin{array}{cc}
M_{\frac{p([N,t])}{4},0}^{0}\#\Sigma _{\frac{p([N,t])}{24}} & \text{if }%
\varepsilon ([N,t])=0 \\ 
M_{\frac{p([N,t])}{4},0}^{0} & \text{if }\varepsilon ([N,t])=1%
\end{array}%
\right. $
\end{quote}

\noindent which completes the proof.$\square $

\bigskip

\noindent \textbf{Proof of Lemma 3.4}. The construction of $W$ and the
corresponding calculations will be divided into two cases depending on $%
\varepsilon ([N,t])=0$ or $1$. Let $\pi _{0}:W_{t}\rightarrow N$ be the
associated disk bundle of the circle bundle $\pi :N_{t}\rightarrow N$.

Case 1\textbf{.} $\varepsilon ([N,t])=1$: Take an embedding $%
f:S^{2}\hookrightarrow $ Interior $W_{t}$ that represents a generator of the
group $H_{2}(W_{t})\cong 
\mathbb{Z}
$. Since $w_{2}(W_{t})=0$(see Example 3.1), the $f$ extends to an embedding $%
\overline{f}:S^{2}\times D^{6}\hookrightarrow $ Interior $W_{t}$. The $8$%
--manifold $W$ is obtained from $W_{t}$ by surgery along $\overline{f}$. On
one hand, it is clear that $\partial W\cong \partial W_{t}\cong N_{t}$. On
the other hand, from the homotopy equivalences

\begin{quote}
$X\simeq W_{t}\cup _{f}D^{3}\simeq W\cup D^{6}$
\end{quote}

\noindent with $X$ the trace of the surgery \cite[P.83-84]{B} we find that $%
W $ is $3$--connected with $H_{4}(W)\cong H_{4}(X)\cong 
\mathbb{Z}
$. Moreover, since $N_{t}$ is $2$--connected and $W$ is $3$--connected we
can conclude that $H_{i}(W)\cong H^{8-i}(W,N_{t})=0$ for $i\geq 5$ by the
Lefschetz duality\ and the cohomology exact sequence. Hence from the
Whitehead theorem we get $W\simeq S^{4}$.

Case 2\textbf{.} $\varepsilon ([N,t])=0$: The desired manifold $W$ is
constructed as follows. Represent the generator $x\cap \lbrack N]\in
H_{2}(N) $ by an embedding $\overline{g}:S^{2}\times D^{4}\hookrightarrow N$
($w_{2}(N)=0$). Let $\overline{h}$ be the pull-back of $\overline{g}$ by the
projection $\pi $ as shown in the diagram

\begin{enumerate}
\item[(3.1)] $%
\begin{array}{ccc}
S^{3}\times D^{4} & \overset{\overline{h}}{\hookrightarrow } & N_{t} \\ 
\downarrow &  & \pi \downarrow \\ 
S^{2}\times D^{4} & \overset{\overline{g}}{\hookrightarrow } & N%
\end{array}%
$,
\end{enumerate}

\noindent and let $\widetilde{W}:=N_{t}\times \lbrack 0,1]\cup _{(\overline{h%
},1)}D^{4}\times D^{4}$. Then

\begin{enumerate}
\item[(3.2)] $\partial \widetilde{W}\cong N_{t}\sqcup ($--$\Sigma _{r})$ for
some $r\in 
\mathbb{Z}
$
\end{enumerate}

\noindent since in (3.1) the map $\overline{h}$ induces an isomorphism $\pi
_{3}(S^{3}\times D^{4})\rightarrow \pi _{3}(N_{t})$ \cite[Lemma 1]{W2}. The
manifold $W$ is obtained from $\widetilde{W}$ by removing the tubular
neighborhood of a smooth arc $\alpha :$ $[0,1]\rightarrow \widetilde{W}$
with $\alpha (0)\in N_{t}$, $\alpha (1)\in \partial W^{^{\prime }}$ and $%
\alpha (0,1)\subset $Interior $\widetilde{W}$.

It remains to show that

\begin{quote}
i) $W\simeq S^{4}$; ii) $r=\frac{p([N,t])}{24}$(in (3.2)).
\end{quote}

\noindent The property i) follows from the facts that the trace $\widetilde{W%
}$ of the surgery along $\overline{h}$ has the homotpy type $\Sigma _{r}\cup
D^{4}$ and the homeomorphism type of $W$ is obtained from $\widetilde{W}$ to
collapse the component $\Sigma _{r}$ of $\partial \widetilde{W}$ to a point.

For the property ii) we only need to show $\mu (\Sigma _{r})\equiv \frac{p(N)%
}{24\cdot 28}\func{mod}%
\mathbb{Z}
$. By the collar neighborhood theorem, the homotopy sphere $\Sigma _{r}$ in
(3.2) bounds an 8--manifold $W^{^{\prime }}:=W_{t}\cup _{\overline{h}%
}D^{4}\times D^{4}$. For the convenience of calculation, we make use of an
alternative decomposition $W^{^{\prime }}=W_{t}\cup _{\overline{i_{0}}}%
\mathbb{C}
P^{2}\times D^{4}$ where $\overline{i_{0}}$ is the pull-back of $\overline{g}
$ by the projection $\pi _{0}$ as in the diagram

\begin{quote}
$%
\begin{array}{ccc}
V\times D^{4} & \overset{\overline{i_{0}}}{\hookrightarrow } & W_{t} \\ 
\downarrow &  & \pi _{0}\downarrow \\ 
S^{2}\times D^{4} & \overset{\overline{g}}{\hookrightarrow } & N%
\end{array}%
$.
\end{quote}

From the isomorphism

\begin{quote}
$%
\begin{array}{ccccc}
H_{4}(%
\mathbb{C}
P^{2}) & \oplus & H_{4}(W_{t}) & \underset{\cong }{\overset{i_{1\ast }\oplus
i_{2\ast }}{\rightarrow }} & H_{4}(W^{^{\prime }})%
\end{array}%
$
\end{quote}

\noindent by the Mayer--Vietoris sequence with $\alpha \in H_{4}(W_{t})\cong 
\mathbb{Z}
$ the generator satisfies $\left\langle \pi _{0}^{\ast }x,\alpha
\right\rangle =1$ (see Section 2) and $i_{1}:%
\mathbb{C}
P^{2}\rightarrow W^{^{\prime }},i_{2}:W_{t}\rightarrow W^{^{\prime }}$ the
inclusions, we can see below that the intersection matrix of $W^{^{\prime }}$
with respect to the basis $x_{1},x_{2}\in H^{4}(W^{^{\prime }},\partial
W^{^{\prime }})$ is

\begin{center}
$\left( 
\begin{array}{cc}
0 & 1 \\ 
1 & 0%
\end{array}%
\right) $
\end{center}

\noindent where $i_{1\ast }[%
\mathbb{C}
P^{2}]=x_{1}\cap \lbrack W^{^{\prime }},\partial W^{^{\prime }}],i_{2\ast
}\alpha =x_{2}\cap \lbrack W^{^{\prime }},\partial W^{^{\prime }}]$. First
observe that

\begin{quote}
$\left\langle x_{1}\cup x_{1},[W^{^{\prime }},\partial W^{^{\prime
}}]\right\rangle =0$
\end{quote}

\noindent as the normal bundle of $i_{1}$ is trivial. Next since the
self-intersection number of $i_{2\ast }\alpha $ is the same as that of $%
\alpha $ and the homomorphism $j^{\ast }:H^{4}(W_{t},N_{t})\rightarrow
H^{4}(W_{t})$ is trivial (see Example 3.1), this implies

\begin{quote}
$\left\langle x_{2}\cup x_{2},[W^{^{\prime }},\partial W^{^{\prime
}}]\right\rangle =\left\langle j^{\ast }D_{W_{t}}\alpha \cup D_{W_{t}}\alpha
,[W_{t},N_{t}]\right\rangle =0$\cite[p.115]{B}
\end{quote}

\noindent where $D_{W_{t}}\alpha $ denotes the Lefschetz duality of $\alpha $%
. Finally we have

\begin{center}
$\left\langle x_{1}\cup x_{2},[W^{^{\prime }},\partial W^{^{\prime
}}]\right\rangle =\left\langle j^{^{\prime }\ast }x_{1},i_{2\ast }\alpha
\right\rangle =\left\langle i_{2}^{\ast }j^{^{\prime }\ast }x_{1},\alpha
\right\rangle =\left\langle \pi _{0}^{\ast }x,\alpha \right\rangle =1$
\end{center}

\noindent with $j^{^{\prime }}:W^{^{\prime }}\rightarrow (W^{^{\prime
}},\partial W^{^{\prime }})$ the inclusion and where $\pi _{0}^{\ast
}x=i_{2}^{\ast }j^{^{\prime }\ast }x_{1}$ follows from the relations $\pi
_{0}^{-1}f[S^{2}]=i_{2}^{-1}i_{1}[%
\mathbb{C}
P^{2}]$, $x=D_{N}f_{\ast }[S^{2}]$ and $x_{1}=D_{W^{^{\prime }}}i_{1\ast }[%
\mathbb{C}
P^{2}]$. Thus the signature $\sigma (W^{^{\prime }})$ is $0$.

We can take $z$ to be a generator of $H^{2}(W^{^{\prime }})\cong 
\mathbb{Z}
$ since $i_{2}^{\ast }w_{2}(W^{^{\prime }})=w_{2}(W_{t})\neq 0$ by Example
3.1 and the isomorphism $i_{2}^{^{\ast }}TW^{^{\prime }}\cong TW_{t}$ . To
get the values of $z^{2},p_{1}(W^{^{\prime }})$, it is necessary to compute
the images of $z^{2},p_{1}(W^{^{\prime }})$ under the isomorphism

\begin{quote}
$i_{1}^{\ast }\oplus i_{2}^{\ast }:H^{4}(W^{^{\prime }})\rightarrow H^{4}(%
\mathbb{C}
P^{2})\oplus H^{4}(W_{t})$,
\end{quote}

\noindent whose matrix with respect to the basis $\{j^{^{\prime \ast
}}x_{1},j^{^{\prime \ast }}x_{2}\}$ and $\{[%
\mathbb{C}
P^{2}]^{\ast },\pi _{0}^{\ast }x\}$ is the same as the intersection matrix
of $W^{^{\prime }}$ with respect to the basis $x_{1},x_{2}$, where $[%
\mathbb{C}
P^{2}]^{\ast }\in H^{4}(%
\mathbb{C}
P^{2})$ satisfies $\left\langle [%
\mathbb{C}
P^{2}]^{\ast },[%
\mathbb{C}
P^{2}]\right\rangle =1$. Since $i_{1}^{\ast }z\in H^{2}(%
\mathbb{C}
P^{2}),i_{2}^{\ast }z\in H^{2}(W_{t})$ are generators, it follows that

\begin{quote}
$i_{1}^{\ast }\oplus i_{2}^{\ast }(z^{2})=(i_{1}^{\ast }z^{2},\pi _{0}^{\ast
}t^{2})=([%
\mathbb{C}
P^{2}]^{\ast },0)$.
\end{quote}

\noindent Moreover, according to the isomorphisms $i_{2}^{^{\ast
}}TW^{^{\prime }}\cong TW_{t}$ and $i_{1}^{^{\ast }}TW^{^{\prime }}\cong T%
\mathbb{C}
P^{2}\oplus \varepsilon ^{4}$, the relations $\left\langle p_{1}(%
\mathbb{C}
P^{2}),[%
\mathbb{C}
P^{2}]\right\rangle =3$ and $p_{1}(W_{t})$ $=p([N,t])\pi _{0}^{\ast }x$
imply that

\begin{quote}
$i_{1}^{\ast }\oplus i_{2}^{^{\ast }}p_{1}(W^{^{\prime }})=(3[%
\mathbb{C}
P^{2}]^{\ast },p([N,t])\pi _{0}^{\ast }x)$.
\end{quote}

\noindent Therefore we can see that

\begin{quote}
$z^{2}=j^{^{\prime \ast }}x_{2};p_{1}(W^{^{\prime }})=p([N,t])j^{^{\prime
\ast }}x_{1}+3j^{^{\prime \ast }}x_{2}$.
\end{quote}

\noindent Again from these relations and the intersection form of $%
W^{^{\prime }}$ we get

\begin{quote}
$p_{1}^{2}(W^{^{\prime }})=6p([N,t]);z^{2}p_{1}(W^{^{\prime
}})=p([N,t]);z^{4}=0.$
\end{quote}

\noindent Consequently, substituting these values in the formula of $\mu ,$
this implies

\begin{quote}
$\mu (M^{^{\prime }})=\frac{p(N)}{24\cdot 28}\func{mod}%
\mathbb{Z}
$.$\square $
\end{quote}

\noindent \textbf{Remark 3.6}. In a communication concerning this work
Diarmuid Crowley pointed out that according to a result of Wilkens \cite[%
Theorem 1 (ii)]{Wl2} the decomposition $N_{t}\cong M_{l,k}^{0}\#\Sigma _{r}$
in Lemma 3.4 can be simplified as $N_{t}\cong M_{l,k}^{0}$ when $k([N,t])=0$%
, which will play a role in the proof of Corollary 4.5 in the coming section.

In the recent paper \cite{CN} (see also \cite{C1}\cite{CN2}) Crowley and
Nordstrom generalised the classical Eells-Kuiper invariant $\mu $. Their new
invariant can be applied to give a simple proof of the diffeomorphism $%
N_{t}\cong M_{l,k}^{0}$ when $k([N,t])=0$.

\section{Proof of Theorem 1.3 and applications}

We establish Theorem 1.3 and present some applications.

\bigskip

\noindent \textbf{Proof of Theorem 1.3. }Let $M$ be a $2$--connected $7$%
--manifold with a regular circle action.\textbf{\ }By Lemmas 2.1 and 2.2 $M$
is the total space of the oriented circle bundle over $N\#_{r}S^{3}\times
S^{3}$ with Euler class $\overline{t}\in H^{2}(N\#_{r}S^{3}\times
S^{3})\cong \mathbb{Z}$ a generator, where $[N,t]\in \Theta $, $r\in 
\mathbb{N}
$. Identify $\overline{t}$ with the generator $t\in H^{2}(N)\cong \mathbb{Z}$
under the isomorphism $H^{2}(N)\rightarrow H^{2}(N\#_{r}S^{3}\times S^{3})$
induced by the map $N\#_{r}S^{3}\times S^{3}\rightarrow N$ collapsing $%
\#_{r}S^{3}\times S^{3}$ to a point. By Lemmas 3.3 and 3.4 it suffices to
show that $M\cong N_{t}\#_{2r}S^{3}\times S^{4}$.

Consider the decomposition

\begin{quote}
$N\#_{r}S^{3}\times S^{3}=(N\backslash \overset{\circ }{D}_{1})\cup
_{f}(\#_{r}S^{3}\times S^{3}\backslash \overset{\circ }{D}_{2})$
\end{quote}

\noindent with $D_{i}\cong D^{6}$ and $f:\partial D_{2}\rightarrow $ $%
\partial D_{1}$ a diffeomorphism. Since the restriction of the bundle $%
N_{t}\rightarrow N$ on $D_{1}$ is trivial and $N_{t}\cong N_{t}\#S^{7}$ one
has the corresponding decomposition

\begin{quote}
$M\cong (N_{t}\backslash \overset{\circ }{D}_{1}\times S^{1})\cup _{f\times
id}((\#_{r}S^{3}\times S^{3}\backslash \overset{\circ }{D}_{2})\times
S^{1})\cong N_{t}\#M_{0}$
\end{quote}

\noindent where $id$ is the identity on $S^{1}$, and where

\begin{quote}
$M_{0}=(S^{7}\backslash \overset{\circ }{D}_{1}\times S^{1})\cup _{f\times
id}((\#_{r}S^{3}\times S^{3}\backslash \overset{\circ }{D}_{2})\times
S^{1}). $
\end{quote}

\noindent Since $M_{0}$ can be easily identified with the total space of the
oriented circle bundle over $CP^{3}\#_{r}S^{3}\times S^{3}$ with Euler class
a proper generator of $H^{2}(CP^{3}\#_{r}S^{3}\times S^{3})\cong \mathbb{Z}$%
, a calculation similar to that in Example 3.1 shows that the invariant
system $\{H,\frac{p_{1}}{2},b,\mu \}$ for $M_{0}$ and $\#_{2r}S^{3}\times
S^{4}$ coincides. Consequently $M_{0}$ is diffeomorphic to $%
\#_{2r}S^{3}\times S^{4}.$ This shows that $M\cong N_{t}\#_{2r}S^{3}\times
S^{4}$ which completes the proof.$\square $

\bigskip

A classical topic is to decide which homotopy spheres admit smooth regular
circle actions (\cite{Hs} \cite{L} \cite{MY} \cite{S1}). Combining Theorem
1.3 with Example 3.2 we recover the classical computation of Montgomery and
Yang \cite{MY} .

\bigskip

\noindent \textbf{Corollary 4.1.} Among the $28$ homotopy $7$--spheres $%
\Sigma _{r},0\leq r\leq 27$ the following ones admit smooth regular circle
actions

\begin{quote}
$\Sigma _{r},r=0,4,6,8,10,14,18,20,22,24$.$\square $
\end{quote}

In term of our notation the unit tangent bundle of the sphere $S^{4}$ is $%
M_{-1,2}^{0}$. The additive property of the Eells-Kuiper invariant $\mu $
shows that $M_{-1,2}^{0}\#\Sigma _{r}$ with $0\leq r\leq 27$ represent all
the diffeomorphism types of the smooth manifolds homeomorphic to $%
M_{-1,2}^{0}$. One can deduce from Theorem 1.3 and Example 3.2 that

\bigskip

\noindent \textbf{Corollary 4.2}. All the smooth manifolds homeomorphic to
the unit tangent bundle of the sphere $S^{4}$ and admitting smooth regular
circle actions are

\begin{quote}
$M_{-1,2}^{0}\#\Sigma _{r},r=0,2,6,7,8,12,14,15,16,19,20,23,26.\square $
\end{quote}

In \cite{GVZ} Grove, Verdiani and Ziller constructed on the manifold $%
M_{-1,2}^{0}\#\Sigma _{27}$ a metric with positive sectional curvature (see
Goette\cite[p.34-35]{Goe}). According to Corollary 4.2 this manifold does
not admit any smooth regular circle action.

\bigskip

\noindent \textbf{Definition 4.3 }Two regular (resp. smooth regular) circle
actions

\begin{center}
$S^{1}\times M_{i}\rightarrow M_{i},i=1,2,$
\end{center}

\noindent on two manifolds (resp. smooth manifolds) $M_{i}$ are called 
\textsl{equivalent} if there is a equivariant homeomorphism (diffeomorphism) 
$f:M_{1}\rightarrow M_{2}$. Let $\rho _{T}(M)$ (resp.$\rho _{S}(M)$) be the
number of all equivalence classes of regular (resp. smooth regular) circle
actions on a given manifold (resp. smooth manifold) $M$.

\bigskip

Since the number $\rho _{T}(M)$ can be seen as the number of those elements $%
[N,t]\in \Theta $ satisfying $N_{t}\cong M$, we get from Lemmas 2.2 and 3.3
that

\bigskip

\noindent \textbf{Corollary 4.4.} For the family

\begin{quote}
$M=M_{6m,(1+c)k}^{c}\#_{2r}S^{3}\times S^{4},c\in \{0,1\},r\in 
\mathbb{N}
,m,k\in 
\mathbb{Z}
$
\end{quote}

\noindent of manifolds that represent all homeomorphism classes of the $2$%
--connected $7$--manifolds with regular circle actions (see Theorem 1.3) we
have

\begin{quote}
$\rho _{T}(M)=\left\{ 
\begin{array}{cc}
1 & \text{if }k=0\text{ and }m\equiv 1\func{mod}2, \\ 
2 & \text{if }k=0\text{ and }m\equiv 0\func{mod}2, \\ 
\infty & \text{if }k\neq 0%
\end{array}%
\right. .\square $
\end{quote}

Similarly, in the smooth category we get from Lemmas 2.2 and 3.5, together
with Remark 3.6, that

\bigskip

\noindent \textbf{Corollary 4.5}. For the family

\begin{center}
$M=M_{6(1+a)m,(1+a)k}^{0}\#\Sigma _{(1-a)m}\#_{2r}S^{3}\times S^{4},a\in
\{0,1\},r\in 
\mathbb{N}
,m,k\in 
\mathbb{Z}
$
\end{center}

\noindent of manifolds that represent all diffeomorphism classes of the
smooth $2$--connected $7$--manifolds with smooth regular circle actions (see
Theorem 1.3) we have

\begin{quote}
$\sigma (M)=\left\{ 
\begin{array}{cc}
1 & \text{if }k=0,a=0\text{ and }m\equiv 1\func{mod}2, \\ 
2 & \text{if }k=0\text{ and }(1+a)m\equiv 0\func{mod}2, \\ 
\infty & \text{if }k\neq 0%
\end{array}%
\right. .\square $\bigskip
\end{quote}

\noindent \textbf{Acknowledgement} The author is grateful to the referee for
many improvements over the first version of this paper. In particular, the
results in Corollaries 4.2, 4.4 and 4.5 are suggested by him.

The author would also like to thank to Haibao Duan for bring her attention
to the topic, and to Diarmuid Crowley for communication concerning this work
(see Remark 3.6). Thanks are also due to Yang Su and Yueshan Xiong for
valuable discussions.

\end{document}